\documentclass{amsart}
\usepackage{amssymb}
\usepackage{amscd,enumerate,amsfonts,calc,amsmath,verbatim,xypic}

\newcommand{\Ext}{\operatorname{Ext}}
\newcommand{\TExt}{\operatorname{E\widehat{\vphantom{E}x}t}\!{\vphantom  E}}
\newcommand{\TTor}{\operatorname{T\widehat{\vphantom{T}o}r}\!{\vphantom T }}
\newcommand{\Tor}{\operatorname{Tor}}
\newcommand{\Hom}{\operatorname{Hom}}
\newcommand{\Soc}{\operatorname{Socle}}
\newcommand{\HH}{\operatorname{H}}
\newcommand{\Ima}{\operatorname{Im}}

\newcommand{\rank}{\operatorname{rank}}
\newcommand{\edim}{\operatorname{edim}}

\newcommand{\Syz}{\operatorname{Syz}}
\newcommand{\cls}{\operatorname{cls\,}}
\newcommand{\ov}{\overline}

\newcommand{\ges}{\geqslant}
\newcommand{\Ker}{\operatorname{Ker}}
\newcommand{\Hilb}{\operatorname{Hilb}}
\newcommand{\Gdim}{\operatorname{G-dim}}

\newtheorem*{Theorem}{Theorem}

\newcommand{\lf}{\otimes^{\mathbf L}}
\newcommand{\rh}{{\mathbf R}\!\Hom}

\newcommand{\fm}{{\mathfrak m}}

\newcommand{\ci}{({\bf ci})}
\newcommand{\gor}{({\bf gor})}
\newcommand{\te}{({\bf te})}
\newcommand{\et}{({\bf et})}
\newcommand{\ee}{({\bf ee})}
\newcommand{\gap}{({\bf gap})}
\newcommand{\ab}{({\bf ab})}
\newcommand{\au}{({\bf ac})}
\newcommand{\ac}{({\bf uac})}
\newcommand{\tv}{({\bf tv})}
\newcommand{\bd}{\boldsymbol}

\theoremstyle{remark}

\swapnumbers

\theoremstyle{plain}
\newtheorem{theorem}{Theorem}[section]

\newtheorem{proposition}[theorem]{Proposition}
\newtheorem{lemma}[theorem]{Lemma}
\newtheorem{corollary}[theorem]{Corollary}

\theoremstyle{definition}
\newtheorem{chunk}[theorem]{}
\newtheorem{subchunk}{}
\newtheorem*{chunk*}{}

\theoremstyle{remark} 
\newtheorem*{Claim}{Claim}

\newtheorem{remark}[theorem]{Remark}
\newtheorem{remarks}[theorem]{Remarks}

\numberwithin{equation}{theorem}

\numberwithin{subchunk}{theorem}

\begin{document}\title[Nonvanishing Cohomology]
{Nonvanishing Cohomology\\
and Classes of Gorenstein rings} 
\author[D.~A.~Jorgensen]{David A.~Jorgensen}
 \address{Department of Mathematics, University of Texas at Arlington, 
Arlington, TX 76019}
\email{djorgens@math.uta.edu}
\author[L.~M.~\c Sega]{Liana M.~\c Sega}
\address{Mathematical Sciences Research Institute, 1000 Centennial Drive,
Berkeley, CA  94720-5070}
\email{lsega@msri.org}

\begin{abstract} We give counterexamples to the following conjecture
  of Auslander: given a finitely generated module $M$ over an Artin
  algebra $\Lambda$, there exists a positive integer $n_M$ such that
  for all finitely generated $\Lambda$-modules $N$, if
  $\Ext_{\Lambda}^i(M,N)=0$ for all $i\gg 0$, then
  $\Ext_{\Lambda}^i(M,N)=0$ for all $i\geq n_M$. Some of our examples
  moreover yield homologically defined classes of commutative local
  rings strictly between the class of local complete intersections 
and the class of local Gorenstein rings.
\end{abstract}

\thanks{This work was done while the first author was a member and the
  second author was a postdoctoral fellow at the Mathematical Sciences
  Research Institute, Berkeley, CA.  They are both grateful to MSRI
  for its generous support.}
\date{\today} 
\maketitle

\section*{introduction} In this paper we give examples on the
vanishing of Ext and Tor which simultaneously disprove a conjecture of
Auslander and identify new classes of commutative local (meaning also
Noetherian) rings lying
between the class of local complete intersections and the class of
local Gorenstein rings.

The following conjecture of Auslander appears in \cite[p.\ 795]{A} and
\cite{H}:
{\it let $\Lambda$ be an Artin algebra.  For every finitely generated
$\Lambda$-module $M$ there exists an integer $n_M$ such that for all
finitely generated $\Lambda$-modules $N$, if
$\Ext^i_\Lambda(M,N)=0$ for all $i\gg 0$, then 
$\Ext^i_\Lambda(M,N)=0$ for all $i\ge n_M$.}

%As described in \cite{A} and \cite{H}, the Finitistic Dimension
%Conjecture holds for a finite dimensional algebra $\Lambda$ over a
%field $k$ if the conjecture above holds for the enveloping algebra
%$\Lambda^e$. 

Auslander's conjecture is known to hold, for example,
when $\Lambda$ is a group ring of a finite group over a field, by
\cite[2.4]{BCR}, or when $\Lambda$ is a local complete intersection
(see the discussion later in the introduction).
Part (1) of our main theorem below gives a counterexample to
Auslander's Conjecture over a commutative selfinjective Koszul
$k$-algebra.  Part (2) is relevant in the context of recent research
on refinements of the Gorenstein condition, as we shall explain
shortly.

\begin{Theorem} Let $k$ be a field which is not algebraic over a
finite field. Then there exist commutative finite dimensional
selfinjective Koszul $k$-algebras $A$ with Hilbert series $\sum
(\rank_kA_i)t^i=1+5t+5t^2+t^3$ and
  finitely generated graded $A$-modules $M$ with linear resolution and
$\sum(\rank_kM_i)t^i=2t+8t^2+2t^3$
such that the following hold:
\begin{enumerate}[\quad\rm(1)]
\item For each positive integer $q$ there exists a finitely generated
graded $A$-module $N_q$ with linear resolution and
$\sum(\rank_k(N_q)_i)t^i=1+t$ satisfying
$$
\Ext^i_{A}(M,N_q)\ne 0 \quad \text{if and only if} \quad i=0,\ q-1,\ q\,. 
$$

\item There exists a finitely generated graded $A$-module $V$
with linear resolution and  
$\sum(\rank_kV_i)t^i=2+2t$ satisfying
$$
\Tor_i^{A}(M,V)=0 \quad\text{for all \quad$i>0$}\quad\text{and}\quad 
\Ext^i_{A}(M,V)\ne 0 \quad\text{for all \quad$i>0$}\,.
$$
\end{enumerate}
\end{Theorem}

It is easy to see that Auslander's Conjecture holds when 
$\Lambda$ is a commutative local
ring with maximal ideal $\fm$ satisfying $\fm^2=0$. 
We show in a corollary of our main theorem (Corollary
3.3.(2)) that the conjecture already fails for a commutative local
ring $(B,\fm)$ with $\fm^3=0$.

The rings $A$ in our main theorem and the rings $B$ in the corollary 
are constructed by Gasharov and Peeva in \cite{GP} to give a 
counterexample to an unrelated conjecture of Eisenbud.  We turned
to these rings since they admit modules of infinite complete intersection
dimension (see \cite{AGP} for the definition), which is a necessary property of any
module yielding a counterexample to Auslander's Conjecture (cf. \cite{AY}).

Recall that a {\it local complete intersection} is a local ring $R$ 
whose completion with respect to the maximal
ideal $\fm$ is a quotient of a regular local ring by a regular
sequence. Let $\mathcal{CI}$ denote the class of all such rings.  Let
$\mathcal{TE}$ represent the class of commutative local rings $R$
which have the following property: $\Tor_i^R(M,N)=0$ for all $i\gg 0$
implies $\Ext^i_R(M,N)=0$ for all $i\gg 0$, for all finitely generated
$R$-modules $M$ and $N$. One of the main theorems of Avramov and
Buchweitz in \cite{AB} gives the first inclusion in the chain
\[
\mathcal{CI}\subseteq\mathcal{TE}\subseteq\mathcal{GOR}\,. 
\]
(The second inclusion is clear: just take $N=R$ in the definition of
the class $\mathcal{TE}$.)  In \cite{AB}, the authors remark that
forty years of research in commutative algebra have not produced a
class of local rings intermediate between $\mathcal{CI}$ and
$\mathcal{GOR}$, and ask whether either inclusion above is strict. In
a recent paper, Huneke and Jorgensen \cite{HJ} prove that the first
inclusion is strict.  Part (2) of our main theorem shows that so is
the second.

In \cite{HJ} an {\it AB ring\/} is defined to be a local Gorenstein
ring $R$ with the property that $\Ext^i_R(M,N)=0$ for all $i\gg 0$
implies $\Ext^i_R(M,N)=0$ for all $i>\dim R$, for all finitely
generated $R$-modules $M$ and $N$.  Let $\mathcal{AB}$ denote the class
of all AB rings. It is shown in \cite{AB} (cf.\ also \cite{HJ}) that
$\mathcal{CI}\subseteq\mathcal{AB}$, and subsequently in \cite{HJ}
(cf.  also \cite{S}) it is shown that this inclusion is strict.  Note
that the condition defining AB rings is a strengthening of Auslander's
Conjecture.  Therefore part (1) of our main theorem also shows that
$\mathcal{AB}$ lies properly between $\mathcal{CI}$ and
$\mathcal{GOR}$.

The paper is organized as follows. In Section \ref{ap} we give some
positive results on Auslander's Conjecture for commutative
non-Gorenstein rings.  In particular, we prove that Golod rings, and
commutative local rings which are ``small'' in various senses satisfy
(a strong form of) Auslander's Conjecture.

The rings $A$ and $B$ and the corresponding modules are defined in
Section \ref{prelims}.  Here we also explain our method for computing
homology and cohomology.

The main theorem above is an immediate consequence of Corollary
\ref{counterexamples}(1) and Proposition \ref{2} proved in Section
\ref{proofs}. We also give there the corresponding statements for the
ring $B$, and compare these with the results of Section \ref{ap},
noting that our examples involving the ring $B$ are ``smallest'' in
various senses where one can expect Auslander's Conjecture to fail.

In the final Section \ref{new rings} we discuss classes of
homologically defined local Gorenstein rings, including the ones
described above.  We give local Gorenstein rings which are known
to satisfy Auslander's Conjecture, and we compare these rings to our
examples from Section \ref{proofs}.

\section{Some commutative rings for which Auslander's Conjecture holds}
\label{ap}

In this section $R$ denotes a commutative local ring, with maximal
ideal $\fm$ and residue field $k$.

As is evidenced by results of \cite{AB} and \cite{HJ}, Auslander's
Conjecture is relevant and interesting in the context of commutative
local rings (of possibly nonzero Krull dimension).
It turns out that all the commutative local rings for which 
Auslander's Conjecture is
known to hold actually satisfy a stronger condition, which we call the
{\it Uniform Auslander Condition\/}:

\medskip

\noindent\ac\quad {\it There exists 
an integer $n$ such that for all finitely generated $R$-modules $M$ and
$N$, if $\Ext^i_R(M,N)=0$ for all $i\gg 0$ then
$\Ext^i_R(M,N)=0$ for all $i\geq n$.}  

\medskip

In this section we prove that $\ac$ holds for certain rings which are
small in various senses. Let $\edim R$ denote the minimal number of
generators of $\fm$ and $\lambda(R)$ denote the length of $R$.

\begin{proposition}
\label{minimal} The local ring $(R,\fm)$ satisfies $\ac$ if  any
one of the following conditions holds.

\begin{enumerate}
[\quad\rm(1)]
\item $\fm^2=0$.
\item $\edim R-\dim R\le 2$.
\item $\fm^3=0$ and $\edim R=3$.
\item $\fm^3=0$ and $\lambda(R)\le 7$.
\end{enumerate} 
\end{proposition}

Recall that the the {\em Poincar\'e series\/} of $M$ over
 $R$ is the formal power series 
 \[
P^R_M(t)=\sum_{i=0}^{\infty}b_i(M)t^i\in\mathbb Z[[t]]
 \]
 where $b_i(M)=\rank_k\Tor^R_i(M,k)$ are the {\em Betti numbers} of
 $M$.
 
 Since some of the results existing in the literature are stated in
 terms of Tor rather than Ext, we remind the reader of the following:

\begin{chunk}
\label{dual}
Assume that $R$ is artinian and let $E$ denote the injective hull of
$k$. For an $R$-module $M$ we set $M^\vee=\Hom_R(M,E)$. By Matlis
duality, for all finitely generated $R$-modules $M$ and $N$ and all $i$
we have: 
$$
\Tor_i^R(M,N^\vee)\cong \Ext^i_R(M,N)^\vee.
$$
\end{chunk}

In some cases one can actually prove a property much stronger than $\ac$. 
We call it {\it trivial vanishing} and it states:  

\medskip

\noindent\tv\quad {\it For any pair $(M,N)$ of finitely generated 
$R$-modules, 
if $\Ext^i_R(M,N)=0$ for all $i\gg 0$, then either
$M$ has finite projective dimension or $N$ has finite injective dimension.}

\begin{proof}[Proof of Proposition {\rm \ref{minimal}}]
Let $M$ and $N$ be finitely generated $R$-modules such that 
$\Ext^i_R(M,N)=0$ for all $i\gg 0$ and assume that $M$ is not free.

(1) The first syzygy $\Syz_1(M)$ in a minimal free resolution of $M$
is annihilated by $\fm$, hence it is a finite sum of copies of
$k$. Since $\Ext^i_R(M,N)=0$ for some $i>1$ implies
$\Ext^{i-1}_R(\Syz_1(M),N)=0$, we conclude that $N$ has finite
injective dimension. The ring therefore satisfies $\tv$, and hence
$\ac$.

(2) By Scheja \cite{Sch}, $R$ is either a complete intersection or
a Golod ring (see \ref{defngolod}).  If it is a complete intersection,
apply \cite[4.7]{AB} (cf. also the last section).  If it is Golod, then
apply Proposition \ref{Golod} below.

(3) If $N$ is injective, then $\Ext^i_R(M,N)=0$ for all $i>0$.
Therefore assume that $N$ is not injective, hence $N^\vee$ is
non-free.  From \ref{dual} we have $\Tor_i^R(M,N^\vee)=0$ for all
$i\gg 0$.  By taking syzygies, we may assume that there exist finitely
generated non-zero $R$-modules $X$ and $Y$ such that $\Tor_i^R(X,Y)=0$
for all $i>0$ and $\fm^2X=\fm^2Y=0$. We conclude from \cite[2.5]{HSV}
that there exist positive integers $u$, $v$ such that $u+v=3$ and
$b_{i+1}(X)=ub_i(X)$, $b_{i+1}(Y)=vb_i(Y)$ for all $i\ge 0$. It
follows that one of the modules $X$ or $Y$ has constant Betti numbers
(because either $u=1$ or $v=1$), hence one of the modules $M$ or
$N^\vee$ has eventually constant Betti numbers. Using a result of
Avramov \cite[1.6]{Av} we conclude that one of the modules
$\Syz_1(M)$, $\Syz_1(N^\vee)$ has a periodic resolution of period
$2$. The hypothesis then implies $\Tor_i^R(M,N^\vee)=0$ for all $i>1$,
hence $\Ext^i_R(M,N)=0$ for all $i>1$.

(4) By (3), we may assume that $\edim R\geq 4$. The ring $R$ then
satisfies the condition in the hypothesis of \cite[3.1]{HSV}, hence,
in view also of \ref{dual}, it satisfies $\tv$.
\end{proof}

\begin{chunk}
\label{defngolod}
Serre proved a coefficientwise inequality
 \[
P_{k}^R(t)\preccurlyeq\frac{(1+t)^{\edim R}}{1-\sum_{j=1}^{\infty}\rank
\HH_{j}(K^R)t^{j+1}}
 \]
of formal power series, where $K^R$ denotes the Koszul complex on
a minimal set of generators of $\fm$.  If equality holds, then $R$
is said to be a {\em Golod ring\/}.
\end{chunk}

\begin{proposition}
\label{Golod}
If $R$ is a Golod ring, then it satisfies $\tv$, and hence $\ac$.
\end{proposition}

For the proof we need some considerations on complexes. We refer to
\cite[Appendix A]{AB} for the basic notions. The Poincar\'e series of
a complex is the extension of the corresponding notion for modules,
cf.\ \cite[\S 7]{AB}, for example.

Vanishing of homology over Golod rings was studied by Jorgensen
\cite[3.1]{J1}. His result was extended in \cite[8.3]{ABS} to
complexes with finite homology as follows:

\begin{chunk}
 \label{Jorgensen}
Let $R$ be a Golod ring and $M$ and $N$ complexes with finite homology.

If $\Tor^R_i(M,N)=0$ for all $i\gg 0$, then $P^R_{M}(t)$ or
$P^R_{N}(t)$ is a Laurent polynomial.
 \end{chunk}

\begin{chunk}A complex of $R$-modules $D$ is said to be {\em dualizing\/}
if it has
finite homology and there is an integer $d$ such that
$\Ext_R^d(k,D)\cong k$ and $\Ext^i_R(k,D)=0$ for $i\ne d$. (By 
\cite[V, 3.4]{Ha}, this definition agrees with the one given by Hartshorne
in \cite{Ha}.)

Any quotient of a local Gorenstein ring has a
dualizing complex.  In particular, a complete local ring has 
a dualizing complex.

\begin{subchunk}
\label{finitedual}
For a complex $G$ we set
$G^{\dagger}=\rh_R(G,D)$. As noted in \cite[V, \S 2]{Ha}, if $G$ has
finite homology, then so does $G^\dagger$.
\end{subchunk}
\end{chunk}

\begin{proof}[Proof of Proposition {\rm \ref{Golod}}] We may assume 
that $R$ is complete, hence it has a dualizing complex $D$. Let $M$ and
$N$ be finitely generated $R$-modules such that 
$\Ext^i_R(M,N)=0$ for all $i\gg 0$. By \cite[V, 2.6(b)]{Ha} we have:
$$
\rh_R(M,N)^\dagger\simeq M\lf_RN^\dagger
$$

By hypothesis, $\rh_R(M,N)$ has finite homology, and by
\ref{finitedual} so does $M\lf_RN^\dagger$. This means that
$\Tor_i^R(M,N^\dagger)=0$ for all $i\gg 0$, hence $P^R_{M}(t)$ or
$P^R_{N^\dagger}(t)$ is a Laurent polynomial, from \ref{Jorgensen}. It
follows that $M$ has finite projective dimension or $N$ has finite
injective dimension (see \cite[V, 2.6(a)]{Ha}).
\end{proof}

\begin{comment}
\begin{chunk} 
  One would like to eliminate the condition $\fm^3=0$ from the
  statement of (3) and (4) of the proposition. We suspect, but we do
  not have a complete proof, that methods similar to those used in
  \cite{S} would prove that all local rings $R$ with $\edim R-\dim
  R=3$ satisfy $\au$, by analyzing the different classes of such
  rings, described in \cite{Av} in terms of their Poincar\'e series.
  All these classes, except for the class {\bf H}$(p,q)$, can easily
  be seen to satisfy the hypothesis of \cite[1.5]{S}, hence they
  satisfy $\tv$. However, we do not know how to settle the case
  {\bf H}$ (p,q)$ if the polynomial $\text{Den}^R(t)$ indicated in the
  chart at page 50 of \cite{Av} has a positive rational root different
  than $1$. (It is possible that this does not happen at all.)
\end{chunk}
\end{comment}

\medskip

Several classes of local Gorenstein rings are
known to satisfy $\ac$. They will be discussed in
Section \ref{new rings}.

\begin{comment}
\medskip

One non-commutative case where we know Auslander's Conjecture holds
is for group rings $kG$ of a finite group $G$ over a field $k$.
The following is proven in \cite{BCR}.

\begin{chunk} \label{BCR}\cite[2.4]{BCR} Given a finite group $G$, there exists
an integer $r$ such that for any commutative ring $R$ of coefficients
and any $RG$-module $M$, if $\hat H^n(G,M)=0$ for $r+1$ consecutive
values of $n$ then $\hat H^n(G,M)=0$ for all $n$ positive and negative.
\end{chunk}

Here $\hat H^n(G,M)$ denotes Tate cohomolgy.  To see that $\ac$ holds
for $kG$ we recall that $\hat H^n(G,M)\cong H^n(G,M)$, the 
ordinary group cohomology, for all $n>0$, and
$$
H^n(kG,\Hom_R(M,N))\cong\Ext^n_{kG}(M,N)
$$
for all $n$ and all left $kG$-modules $M$ and $N$ \cite[IX, 4.4]{CE}.
Now apply \ref{BCR}.
\end{comment}

\section{Constructions} 
\label{prelims}

Let $k$ be a field and let $\alpha\in k$ be a nonzero element.  In
this section we describe the rings $A_\alpha$ and $B_\alpha$
constructed in \cite{GP} by Gasharov and Peeva, we define the modules
$M$ and $L$ of our main theorem and its corollary, and we discuss our
method for computing homology and cohomology.
\medskip

\noindent{\bf The ring $A=A_{\alpha}$.}

Consider the polynomial ring $k[X_1, X_2, X_3, X_4,X_5]$ in five
(commuting) variables and set $A=k[X_1, X_2, X_3, X_4,X_5]/I_{\alpha}$, 
where $I_{\alpha}$ is the ideal generated by the following quadric
relations:
\begin{gather*}
\alpha X_1X_3+X_2X_3,\; X_1X_4+X_2X_4,\; X_3^2+\alpha X_1 X_5-X_2X_5, \\
X_4^2+X_1X_5-X_2X_5,\; X_1^2,\; X_2^2,\; X_3X_4,\; X_3X_5,\; X_4X_5,\;
X_5^2 \,.
\end{gather*}
By \cite{GP}, $A_{\alpha}$ is a local Gorenstein ring with Hilbert series
$\Hilb_{A_{\alpha}}(t)=1+5t+5t^2+t^3$. As a $k$-vector space, 
it has a basis consisting of the $12$ elements
\[
1,\; x_1,\;x_2,\; x_3,\; x_4,\; x_5,\; x_1x_2,\; x_1x_3,\;
x_1x_4,\; x_1x_5,\; x_2x_5,\; x_1x_2x_5\,,
\]
where $x_i$ denotes the residue class of $X_i$ modulo $I_{\alpha}$.

\medskip

\noindent{\bf The ring $B=B_{\alpha}$.} 

Set $B_\alpha=A_\alpha/(x_5)$. As noted in \cite{GP}, $B_\alpha$ is a
local ring with Hilbert series $\Hilb_{B_\alpha}(t)=1+4t+3t^2$.  As a
$k$-vector space, it has a basis formed by the images in $B_\alpha$ of
the following $8$ elements in $A_\alpha$:
\[
1,\, x_1,\, x_2,\, x_3\, x_4,\, x_1x_2,\, x_1x_3,\, x_1x_4\,.
\]

\medskip

When there is no danger of confusion we will suppress 
$\alpha$ from the notation and simply write $A$ or $B$ for $A_\alpha$
or $B_\alpha$.

\begin{chunk} One may check that the set of generators of $I_{\alpha}$
listed above is itself a Gr\"obner basis for $I_{\alpha}$ with respect
to the reverse lexicographic term order.
Since all of these generators are quadrics, by
\cite[Section 4]{F} we have that $A=A_{\alpha}$ is Koszul.  Similarly,
the generators $$
\alpha x_1x_3+x_2x_3,\ x_1x_4+x_2x_4,\ x_3^2,\ x_4^2,\ x_1^2,\ x_2^2,\ x_3x_4
$$ 
of the ideal defining $B_{\alpha}$ form a Gr\"obner basis of their
ideal with respect to the reverse lexicographic term order, 
and so $B=B_{\alpha}$ is also Koszul.
\end{chunk}

Modules with nonperiodic (or periodic of period $\neq 2$) minimal
resolutions having constant Betti numbers equal to $2$ were given in
\cite{GP} over the rings $A_\alpha$ and $B_\alpha$ with $\alpha\ne \pm
1$. We wanted a module with nonperiodic resolution and constant Betti
numbers, but we found that the modules in \cite{GP} did not provide
counterexamples using our technique.

\medskip

\noindent{\bf The modules $M$ and $L$.}

Let $M$ be the image of the map $d_0\colon A^2\to A^2$ given in the
standard basis of $A^2$ as a free $A$-module by the matrix
\[
\left(\begin{matrix}
x_1& x_3\\
x_4&x_2\end{matrix}\right).
\]
Set $L=M\otimes_AB$. 

\medskip

For any ring $R$ we let $-^*$ denote the $R$-module $\Hom_R(-,R)$.

\begin{lemma}
\label{complete}
Consider the sequence of homomorphisms
$$ 
\bd C:\quad \cdots\to
A^2\xrightarrow{d_{i+1}}A^2\xrightarrow{d_i}A^2\xrightarrow{d_{i-1}}
A^2\to\cdots, 
$$
where each map $d_i$ is given in the standard basis of
$A^2$ over $A$ by the matrix
\[
\left(\begin{matrix}x_1& \alpha^ix_3\\x_4 &x_2\end{matrix}\right)\,.
\]
Then $\bd C$ is an exact complex. Moreover, the complexes $\bd C^*$,
$\bd C\otimes_AB$, and $\bd C^*\otimes_AB$ are also exact.
\end{lemma}

We use the subscript of the
differential to keep track of degrees within our complexes.
 
\begin{remark} 
\label{complete resolution}
Let $W$ be a finitely generated module over a Noetherian ring $R$. A
{\it complete resolution} of $W$ is a complex $\bd T$ of finite
projective $R$-modules such that $\HH_i(\bd T)=\HH_i(\bd T^*)=0$ for
all $i\in \mathbb Z$, and $\bd T_{\ges r}=\bd P_{\ges r}$ for some
projective resolution $\bd P$ of $W$ and some $r$.

If $W\ne 0$, then its {\it G-dimension} is the shortest length of a
resolution by modules $G$ with $G\cong G^{**}$ and
$\Ext^i(G,R)=\Ext^i_R(G^*,R)=0$ for all $i>0$; it is denoted
$\Gdim_RW$. By a basic result of Auslander and Bridger \cite{ABr}, the
ring $R$ is Gorenstein if and only if every finite $R$-module $W$ has
$\Gdim_RW<\infty$.

By \cite[4.4.4]{AB}, $W$ has a complete resolution if and only if
$\Gdim_RW<\infty$. 
Lemma \ref{complete} shows that $\bd C$ is a complete resolution of
the $A$-module $M$ and $\bd C\otimes_AB$ is a complete resolution of
the $B$-module $L$. In particular, $L$ has finite
G-dimension. (Since G-dimension is bounded by the dimension of the
ring, we actually have $\Gdim_BL=0$.)
\end{remark}
 
\begin{proof}[Proof of Lemma {\rm \ref{complete}}] It is immediate to 
check from the defining equations of $A$ that $d_id_{i+1}=0$.

As a $k$-vector space, $A^2$ is 24-dimensional. We let $(a,b)$ denote
an element of $A^2$ written in the standard basis of $A^2$ as a free
$A$-module. It is easy to see that for each $i$, the following
elements in $\Ima d_i$ are linearly independent over $k$:
\begin{align*}
  d_i(1,0)&=(x_1,x_4)     &d_i(x_5,0)&=(x_1x_5,0)\\
  d_i(0,1)&=(\alpha^ix_3,x_2)      &d_i(0,x_1)&=(\alpha^ix_1x_3,x_1x_2)\\
  d_i(x_1,0)&=(0,x_1x_4)    &d_i(0,x_3)&=(0,-\alpha x_1x_3) \\
  d_i(x_2,0)&=(x_1x_2,-x_1x_4)        &d_i(0,x_5)&=(0,x_2x_5)\\
  d_i(x_3,0)&=(x_1x_3,0)     & d_i(x_2x_5,0)&=(x_1x_2x_5,0)         \\
  d_i(x_4,0)&=(x_1x_4, x_2x_5-x_1x_5)& d_i(0,x_1x_5)&=(0,x_1x_2x_5)
\end{align*}
As a consequence, $\rank_k(\Ima d_i)\ge 12$ for each $i$. On the other hand, we
have
$$
\rank_k(\Ker d_{i})=\rank_k(A^2)-\rank_k(\Ima d_{i})\le 12 $$ for each
$i$. It follows that $\rank_k(\Ima d_{i+1})=\rank_k(\Ker d_i)=12$ for
each $i$, hence the complex $\bd C$ is exact.  Since $A$ is
selfinjective, the complex $\bd C^*$ is exact as well.

Similar computations over the ring $B$ show that the complexes $\bd
C\otimes_AB$ and $C^*\otimes_AB$ are exact.
\end{proof}

\begin{chunk} The proof of the lemma shows that
$$
\Hilb_M(t)=\Hilb_{M^*}(t)=2t+8t^2+2t^3\,,
$$
and one may verify that 
$$
\Hilb_L(t)=\Hilb_{L^*}(t)=2t+6t^2\,.
$$
\end{chunk}

Every finitely generated graded module $W$ over a standard graded
local ring $R$ with $R_0=k$ has a minimal graded free resolution. 
Consequently, the modules $\Tor_i^R(W,k)$ inherit a structure of graded
$R$-modules. The
bigraded Poincar\'e series of $W$ is the formal power series in two
variables: 
$$ 
P^R_W(t,z)=\sum_{i,j}\rank_k\Tor_i^R(W,k)_j\,t^iz^j\,, 
$$
where $j$ is the index of degree.
The usual Poincar\'e series is obtained by letting $z=1$.

The module $W$ is said to have a {\em linear resolution} if all its
minimal generators are in the same degree $p$ and $W$ has a minimal
graded free resolution in which all the entries of the matrices
defining the differentials have degree $1$. This is equivalent to 
$\Tor_i(W,k)_j=0$ for all $j\ne i+p$.

By definition, the $k$-algebra $R$ is Koszul if
the $R$-module $k$ has a linear resolution. In this case, it is known that
$$
P^R_k(t,z)=\frac{1}{\Hilb_R(-tz)}\,.
$$

\begin{chunk} 
\label{poincare}
It is clear from Lemma \ref{complete} that the $A$-modules $M$ and
$M^*$, as well as the $B$-modules $N$ and $N^*$, have a linear
resolution, and $$
P^A_M(t,z)=P^A_{M^*}(t,z)=P^B_L(t,z)=P^B_{L^*}(t,z)=\frac{2z}{1-tz}\,.
$$ 
Note that all the syzygies of these modules have the same
Poincar\'e series.  Since the rings $A$ and $B$ are Koszul, the
Poincar\'e series of $k$ over $A$ and $B$ are $$
P_k^A(t,z)=\frac{1}{1-5tz+5t^2z^2-t^3z^3}\quad\text{and}\quad
P_k^B(t,z)=\frac{1}{1-4tz+3t^2z^2}\,.  $$
\end{chunk}

\medskip

Next we describe our approach to calculating (co)homology over the
rings $A$, respectively $B$, when one of the modules is $M$,
respectively $L$.

\begin{chunk}
\label{computing}
{\bf Computing Ext and Tor.}
  
We set $\bd F=\bd C_{\ge 0}$ and $\bd G=\bd C_{<0}$. Lemma
\ref{complete} shows that $\bd F$ is a minimal free resolution of $M$
over $A$ and $\bd F\otimes_AB$ is a minimal free resolution of $L$
over $B$. Also, ${\bd G}^*$ is a minimal free resolution of $M^*$ over
$A$ and $\bd G^*\otimes_AB$ is a minimal free resolution of $L^*$ over
$B$.

Let $N$ be a finitely generated $A$-module and let $P$ be a finitely
generated $B$-module.

\begin{subchunk}\label{Tor}
The module $\Tor_i^A(M,N)$ is the $i$th homology of the complex 
\[
\bd F\otimes_AN: \quad \dots \to N^2\xrightarrow{ d_i\otimes_AN}
N^2\xrightarrow{
d_{i-1}\otimes_AN}N^2\to\dots \to N^2\xrightarrow{
d_1\otimes_AN} N^2\to 0\,,
\]
that is,
$$
\Tor_i^A(M,N)=\Ker(d_i\otimes_AN)/\Ima(d_{i+1}\otimes_AN)\,.
$$

Similarly, $\Tor_i^B(L,P)$ is the $i$th homology of the complex
$(\bd F\otimes_AB)\otimes_BP=\bd F\otimes_AP$.
\end{subchunk}

\begin{subchunk}\label{Ext(M,N)}
The module $\Ext^i_A(M,N)$ is the $(-i)$th homology of the complex
$\Hom_A(\bd F,N)$ and this complex can be identified with
\[
\bd F^*\otimes_AN: \quad  
0\to N^2\xrightarrow{ d_1^*\otimes_AN}N^2\to\dots \to N^2
\xrightarrow{ d_i^*\otimes_AN}N^2\xrightarrow
{d_{i+1}^*\otimes_AN} N^2\to\cdots
\]
where $\bd F^*=\Hom_A(\bd F,A)$, and for each $i$ the map $d_i^*$ is given in
the standard basis of $A^2$ by the transpose of the matrix
corresponding to $d_i$. Thus 
$$
\Ext^i_A(M,N)=\Ker(d^*_{i+1}\otimes_AN)/\Ima(d^*_i\otimes_AN)\,.
$$ 

Similarly, $\Ext^i_B(L,P)$ is the $(-i)$th homology of the complex
$\bd F^*\otimes_AP$.
\end{subchunk}

\begin{subchunk}\label{Ext(N,M)}
The module $\Ext^i_A(N,M)$ is isomorphic to $\Tor_i^A(M^*,N)^*$, 
and $\Tor_i^A(M^*,N)$ is the $i$th homology of the complex 
\[
\bd G^*\otimes_AN: \quad \dots \to N^2\xrightarrow{ d_{-i}^*\otimes_AN}N^2
\xrightarrow{
d_{-i+1}^*\otimes_AN}N^2\to\dots \to N^2\xrightarrow{
d_{-1}^*\otimes_AN} N^2\to 0\,.
\]
Thus 
$$
\Ext^i_A(N,M)=\left(
\Ker(d^*_{-i}\otimes_AN)/\Ima(d^*_{-i-1}\otimes_AN)\right)^*\,.
$$

Similarly, $\Tor_i^B(L^*,P)$ is the $i$th homology of the complex $\bd
G^*\otimes_AP$.
\end{subchunk}
\end{chunk}

\section{Results on vanishing}
\label{proofs}

In this section we prove the main results stated in the
introduction. Our method and constructions were inspired
by the paper \cite{He} of Heitmann.

We fix $\alpha\in k$ and use the notation introduced in Section
\ref{prelims}.

For each integer $q$ we set $T_q=A/J_q$, where $J_q$ is the ideal of
$A=A_\alpha$ generated by the linear relations
$$x_1-x_2,\; x_1-\alpha^qx_3,\; x_1-x_4,\;x_5\,.$$ Note that $T_q$ is
also a $B$-module.

We let $o(\alpha)$ denote the order of $\alpha$ in the group of units
of $k$, and set 

$$
s=\begin{cases}
               0\qquad\,\,\text{if $o(\alpha)=\infty$}\\ 
               o(\alpha)\quad\text{otherwise.}
    \end{cases}
$$

\begin{proposition}
\label{1}
The following hold for every integer $q$ and every $i>0$:
\begin{enumerate}[\quad\rm(a)]
\item $\Tor_i^A(M,T_q)\ne 0$  if and only if  $i\equiv q-1,q\mod (s)$.
\item $\Ext^i_A(M,T_q)\ne 0$  if and only if  $i\equiv q-1,q\mod (s)$.
\item $\Ext^i_A(T_q,M)\ne 0$  if and only if  
$i\equiv -q,-q-1\mod (s)$. 
\end{enumerate}
\end{proposition}

\begin{proposition}
\label{1'}
The following hold for every integer $q$ and every $i>0$:
\begin{enumerate}[\quad\rm(a)]
\item $\Tor_i^B(L,T_q)\ne 0$  if and only if 
$i\equiv q-1,q \mod (s)$.
\item $\Ext^i_B(L,T_q)\ne 0$  if and only if 
$i\equiv q-1,q\mod (s)$.
\item $\Tor_i^B(L^*, T_q)\ne 0$  if and only if 
$i\equiv -q,-q-1\mod (s)$.
\end{enumerate}
\end{proposition}

\begin{corollary}
\label{counterexamples}
  If $o(\alpha)=\infty$, then the following hold for any integer
  $q>0$:
\begin{enumerate}[\quad\rm(1)]
\item $\Ext^i_{A}(M,T_{q})\ne 0$  if and only if  $i=0,\ q-1,\ q$.
\item $\Ext^i_{B}(L,T_{q})\ne 0$  if and only if  $i=0,\ q-1,\ q$.
\end{enumerate} 
\end{corollary}

\begin{remark}
  The corollary shows that, when $o(\alpha)=\infty$, the rings
  $A=A_{\alpha}$ and $B=B_{\alpha}$ provide counterexamples to
  Auslander's Conjecture.  In view also of \ref{linear} below, the first
  part gives part (1) of the
main theorem in the introduction. 
  \end{remark}

\begin{chunk}
Let $\fm$ denote the maximal
ideal of $B$. The expression
for $\Hilb_B(t)$ given above indicates 
that $\fm^3=0$, $\edim B=4$ and $\lambda(B)=8$.
Comparing these numerical data with the results stated in Proposition
\ref{minimal}, we see that our examples are
minimal primarily with respect to the invariant $\sup\{n\mid \fm^n=0\}$
and secondarily with respect to the invariants $\edim$ and length.
\end{chunk}

\begin{remark}
\label{Tate}
Let $R$ be a Noetherian ring and let $W$ and $N$ be finitely generated
$R$-modules. Assume that $W$ has a complete resolution $\bd T$, as
defined in \ref{complete resolution}. For each $i$ the {\it Tate
  (co)homology} groups are defined by
$$
\TExt^i_R(W,N)=\HH_{-i}\Hom(\bd T, N)\quad\text{and}
\quad\TTor_i^R(W,N)=\HH_i(\bd T\otimes_RN)\,.
$$

If $r$ is as in \ref{complete resolution}, then it is clear that for
all $i>r$ one has
$$
\TExt^i_R(W,N)\cong\Ext^i_R(W,N)\quad\text{and} 
\quad\TTor_i^R(W,N)\cong \Tor_i^R(W,N)\,.
$$
Also, when $\Gdim_RW=0$, one has
$$
\TExt^{-i-1}_R(W,N)\cong \TTor_i^R(W^*,N)\,.
$$
In terms of Tate (co)homology, the propositions above can be
formulated as follows. Let $q$ be an integer.  Then for all $i$ we have:
\begin{enumerate}[\quad\rm(1)]
\item $\TExt^i_A(M, T_q)=0$ if and only if $i\equiv q-1,q\mod (s)$. 
\item $\TTor_i^A(M,T_q)=0$  if and only if  $i\equiv q-1,q\mod(s)$.
\item $\TExt^i_B(L, T_q)=0$  if and only if  $i\equiv q-1,q\mod (s)$. 
\item $\TTor_i^B(L,T_q)=0$  if and only if  $i\equiv q-1,q\mod(s)$.
\end{enumerate}
\end{remark}

\begin{remark}
  It is now clear that, using the modules in the propositions, one can
  obtain a wide variety of distributions of nonzero (co)homology.  In
  particular, when $s=0$ one can construct arbitrarily large intervals
  of either vanishing or nonvanishing (co)homology: for
  integers $a$, $b$ satisfying $0\leq a<b$ there exist finitely generated
  $A$-modules $N_{a,b}$ and $Z_{a,b}$ such that
\begin{enumerate}[\quad\rm(1)]
\item $\Ext^i_A(M,N_{a,b})\ne 0$ if and only if $a\leq i\leq b$ (and $i=0$), 
and
\item $\Ext^i_A(M,Z_{a,b})=0$ for all $a<i<b$ and
      $\Ext^i_A(M,Z_{a,b})\ne 0$ for $i=a, b$ (and $i=0$). 
\end{enumerate} Indeed, one may take
$N_{a,b}=\bigoplus_{q=a+1}^{b}T_q$ and $Z_{a,b}=T_{a}\bigoplus
T_{b+1}$.

When $s$ is positive, we obtain recurring intervals of
vanishing/nonvanishing cohomology. For example, assume that $s=4$ and
set $T=T_0$.  We have then
\begin{align*}
\Ext^i_A(M,T) &=0 \quad\text{for all}\quad i\equiv 1,2\mod (4) \\
\Ext^i_A(M,T) &\ne 0\quad\text{for all}\quad i\equiv 0,3\mod (4) 
\end{align*}
\end{remark}

\medskip

We give only the proof of Proposition \ref{1}. The proof of
Proposition \ref{1'} is similar.

\begin{proof}[Proof of Proposition {\rm \ref{1}}]
  We let overbars denote residue classes modulo $J_q$ and we perform
  computations of Ext and Tor as explained in
  \ref{computing}.  We only give the proof of (1); the other
  arguments are similar.
 
(1) The differential $\ov d_i=d_i\otimes_AT_q$ of the complex
$\bd F\otimes_AT_q$ is given in the standard basis of $T_q^2$ over $T_q$
by the matrix
\[
\left(\begin{matrix}
\ov x_1& \alpha^{i-q}\ov x_1\\
\ov x_1&\ov x_1\end{matrix}\right).
\]
As a $k$-vector space, $T_q$ has a basis consisting of $1,\ov x_1$
and for each $i\ge 0$
$$
\rank_k(\Ima \ov d_i)=\begin{cases} 1\quad \text{if}\quad i\equiv q\mod(s)\\
2\quad \text{otherwise}.
\end{cases}
$$
Since $\dim_k(\Ker \ov d_i)=\dim_k(T_q^2)-\dim_k(\Ima \ov d_i)$, we then
have
$$
\rank_k(\Ker \ov d_i)=\begin{cases} 3\quad \text{if}\quad i\equiv q\mod(s)\\
2\quad \text{otherwise}.
\end{cases}
$$ 

By \ref{Tor} we have $\Tor_i^A(M,T_q)=\HH_i(\bd F\otimes_AT_q)$
and the conclusion follows from the above computations.
\end{proof}

\begin{chunk} 
\label{linear}
For each $q$ the module $T_q$ has Hilbert series 
$$
\Hilb_{T_q}(t)=1+t.
$$ 

Assume that $o(\alpha)=\infty$. By Proposition
\ref{1} and \ref{poincare} there exists an $A$-module $W$ with
$P^A_W(t,z)=2z(1-tz)^{-1}$ such that $\Tor_i^A(W,T_q)=0$ for all $i>0$
and $W\otimes_RT_q$ is isomorphic to a sum of $2$ copies of $k$ each
generated in degree 1. Indeed, if $q\le 0$ then take $W$ to be a first
syzyzgy of $M$, and if $q>0$ then take $W$ to be a first syzygy of
$M^*$, and use for example \cite[1.4(2)]{HSV} to see that
$W\otimes_RT_q$ is annihilated by the maximal ideal of $A$.

The bigraded version of a usual computation, cf. \cite[1.1]{CM} for
example, gives
\[
P^A_{W\otimes_AT_q}(t,z)=P^A_W(t,z)P^A_{T_q}(t,z)\,. 
\] 
Since $W\otimes_RT_q$ is a sum of copies of $k$, we may use the formula
for $P^A_k(t,z)$ given in \ref{poincare} to conclude that 
\[
P^A_{T_q}(t,z)=\frac{1}{1-4tz+t^2z^2}\,.
\]
The expansion of this fraction as a power series shows that the
$A$-module $T_q$ has a linear resolution.  Similar computations show
that $T_q$ has a linear resolution as a $B$-module as well.
\end{chunk}

\medskip

Now let $U$ be the the cokernel of the map $A^6\to A^2$ given in the
standard bases of $A^6$, respectively $A^2$, over $A$ by the matrix
\[
\left(\begin{matrix}
x_3 &0 &x_1  &x_4 &x_2 &0\\
-x_2&x_3 &-x_4 &0 &0 &x_1
\end{matrix}\right),
\]
and set $V=U\otimes_A A/(\fm^2,x_5)$. Note that $V$ is also a
$B$-module.

The next proposition gives part (2) of the main theorem in
the introduction. (See also \ref{linear2} below). 

\begin{proposition}
\label{2}
With the notation above, the following hold:
\begin{enumerate}[\quad\rm(a)]
\item $\Tor_i^A(M,V)=0$ for all $i>0$.
\item $\Ext^i_A(M,V)\ne 0$ for all $i>0$.
\end{enumerate}
\end{proposition}

\begin{remarks} We can obtain an example of vanishing Exts and
  nonvanishing Tors just by replacing $V$ with $V^*$.  (Since $A$ is
  zero-dimensional and Gorenstein, one has $-^*=-^\vee$. Recall
  then from \ref{dual} that $\Tor_i^A(M,V^*)\cong \Ext^i_A(M,V)^*$ and
  $\Ext^i_A(M,V^*)\cong \Tor_i^A(M,V)^*$.)
  
  One can replace $A$ with $B$ and $M$ with $L$ in the statement of
  Theorem \ref{2}.  The proof is similar.
\end{remarks}

\begin{proof}[Proof of Proposition {\rm \ref{2}}] 
(1). Recall from \ref{Tor} that $\Tor_i^A(M,V)$ is the $i$th
homology of the complex $\bd F\otimes_AV$, where the differential 
$d_i\colon A^2\to A^2$ of $\bd F$ is given in the standard basis of
$A^2$ by the matrix
\[
\left(
\begin{matrix}
x_1& \alpha^ix_3\\ x_4& x_2
\end{matrix}
\right).
\]

The module $V$ is the quotient of $A^2$ by $(\fm^2,x_5)A^2$ and the
relations $(x_3,-x_2)$, $(0,x_3)$, $(x_1,-x_4)$, $(x_4,0)$, $(x_2,0)$,
$(0,x_1)$. We let $\cls a$ denote the residue class in $V$ of $a\in
A^2$. As a $k$-vector space, $V$ has a basis formed by the $4$
elements
$$
v_1:=\cls (1,0),\; v_2:=\cls (0,1),\; v_3:=\cls (x_1,0),\; v_4:=\cls (0,x_2).
$$
Given
this basis for $V$, we then make the obvious choice for a basis
of $V^2$.  

We set $\delta_i=d_i\otimes_AV\colon V^2\to V^2$. For each
$i$, the following elements in $\Ima (\delta_i)$ can be easily
seen to be linearly independent over $k$:
\begin{align*} 
\delta_i(v_1,0)&=(v_3,0)\\
\delta_i(v_2,0)&=(0,v_3)\\
\delta_i(0,v_1)&=\alpha^i(v_4,0)\\
\delta_i(0,v_2)&=(0,v_4).
\end{align*} 
As a consequence, $\rank_k(\Ima \delta_i)\ge 4$. On the other hand,
$$\rank_k(\Ker \delta_i)=\rank_k(V^2)-\rank_k(\Ima \delta_i)\le
8-4=4$$ It follows that $\Ker(\delta_i)=\Ima(\delta_{i+1})$ for each
$i>0$, hence $\HH_i(\bd F\otimes_AV)=0$.

(2). Recall from \ref{Ext(M,N)} that $\Ext_i^A(M,V)$ is the $(-i)$th
homology of the complex $\bd F^*\otimes_AV$, where the differential
$d_i^*\colon A^2\to A^2$ of $\bd F^*$ is given in the standard basis
of $A^2$ by the matrix
\[
\left(
\begin{matrix}
x_1& x_4\\ \alpha^ix_3& x_2
\end{matrix}
\right)
\]
Set $\delta_i^*=d_i^*\otimes_A V$. Similar computations as above show
that only two elements of the basis of $V^2$ over $k$ are not in $\Ker
(\delta_i^*)$. Their images are
\begin{align*}
\delta_i^*(v_1,0)&=(v_3,\alpha^iv_4)\\
\delta_i^*(0,v_2)&=(v_3,v_4)
\end{align*} 
We conclude 
\[
\rank_k(\Ima \delta_i)\le 2\quad \text{and}\quad \rank_k(\Ker
\delta_i)\ge 8-2=6\quad\text{for all $i$.}
\]
We thus have $\HH_i(\bd F^*\otimes_AV)\ne 0$ for all $i>0$.
\end{proof}

\begin{remark} Formulated in terms of Tate (co)homology, the proof of
Proposition \ref{2} shows that
\begin{enumerate}[\quad\rm(a)]
\item $\TTor_i^R(M,V)=0$ for all $i$. 
\item $\TExt^i_R(M,V)\ne 0$ for all $i$. 
\end{enumerate}
\end{remark}

\begin{chunk}
\label{linear2}
The proof of the proposition shows that the module $V$ has 
$$
\Hilb_V(t)=2+2t.
$$
Similar computations as in \ref{linear} show that $V$ has a linear
resolution, both as an $A$-module, 
and as a $B$-module. 
\end{chunk}

For each finitely generated $A$-module $N$ we set 
$$
c(N)=\rank_k(N)-\rank_k(\Soc N).
$$

The next proposition shows that when $o(\alpha)=\infty$, $M$ is 
rather rigid.

\begin{proposition}
\label{3}
When $o(\alpha)=\infty$ the following hold for any finitely generated 
$A$-module $N$:
\begin{enumerate}[\quad\rm(1)]
\item If $\Tor_j^A(M,N)=0$ for some $j>0$, then $\Tor_i^A(M,N)\ne 0$ for
      at most $2c(N)$ values of $i>0$. 
\item If $\Ext^j_A(M,N)=0$ for some $j>0$, then $\Ext^i_A(M,N)\ne 0$ for
      at most $2c(N)$ values of $i>0$. 
\item If $\Ext^j_A(N,M)=0$ for some $j>0$ then $\Ext^i_A(N,M)\ne 0$ for
      at most $2c(N)$ values of $i>0$. 
\item $\Ext^i_A(M,N)=0$ for all $i\gg 0$ if and only if
      $\Ext^i_A(N,M)=0$ for all $i\gg 0$.
\end{enumerate}
\end{proposition}

\begin{remark}
  The statements (1) and (2) remain valid, with similar proofs, when
  replacing $A$ by $B$, and $M$ by $L$.
\end{remark}

\begin{proof}
  
  We will show that if $\HH_j(\bd C^*\otimes_AN)=0$ for some $j$, then
  $\HH_i(\bd C^*\otimes_AN)\ne 0$ for at most $2c(N)$ values of $i$.
  In view of \ref{Ext(M,N)} and \ref{Ext(N,M)}, this proves (2), (3)
  and (4). The proof of (1) is along similar lines, using the complex
  $\bd F\otimes_AN$ and \ref{Tor}, and it is omitted.

For every $i$ set 
\[
u_i=\rank_k\big(\Ima(d^*_{i+1}\otimes_A N)\big)\quad
\text{and}\quad v_i=\rank_k\big(\Ker (d^*_i\otimes_A N)\big)\,.
\]
Since $\bd C^*\otimes_AN$ is a complex, we have $u_i\le v_i$ for all $i$,
with equality if and only if that $\HH_i(\bd C^*\otimes_AN)=0$.

Assume that $u_j=v_j$ for some $j\in\mathbb Z$. We need to show that
$u_i\ne v_i$ for at most $2c(N)$ values of $i\in\mathbb Z$. Set
$$
r=\max\{u_i\mid i\in \mathbb Z\}.
$$
Since $u_i+v_{i+1}=2\rank_kN$ and $u_i\le v_i$ for all $i$, we
conclude that $u_i+u_{i-1}\le 2\rank_k N$, with equality if and only
if $u_i=v_i$. Taking $i=j$ we obtain that $u_j+u_{j-1}=2\rank_kN$.
Since $u_i\le r$ for all $i$, we have $\rank_kN\le r$.

\begin{Claim}
  $u_i\ne r$ for at most $c(N)$ values of $i\in\mathbb Z$.
\end{Claim}

Assuming the claim for the moment, we finish the proof. 

As noted above, we need to show
that $u_i+u_{i-1}\ne 2\rank_kN$ for at most $2c(N)$ values of 
$i\in\mathbb Z$. We have
\[
\{i\in\mathbb Z\mid u_i+u_{i-1}\ne 2\rank_kN\}=
\{i\in\mathbb Z\mid u_i\ne r\}\cup\{i\in\mathbb Z\mid u_{i-1}\ne r\}\,,
\]
and the claim shows that these latter sets have at most $2c(N)$ elements.

\medskip

\noindent {\it Proof of claim.} 
Let $\mathcal E$ be a basis of $N$ over $k$ and set $e=\rank_kN$.  Let
$\chi_i$ be the $e\times e$ matrix which represents in the basis
$\mathcal E$ the map $N\to N$ given by multiplication by $x_i$, for
$i=1,\dots,4$.

With the obvious choice for the basis of $N^2$ over $k$, the map
$d_i\otimes_AN\colon N^2\to N^2$ is represented by the $2e\times 2e$
matrix
\[
\Omega_i=\left(
\begin{matrix}
\chi_1 & \alpha^i\chi_3\\
\chi_4 & \chi_2
\end{matrix}
\right).
\]

Let $k[y]$ be the polynomial ring over $k$ in a single variable $y$.
We consider the following $2e\times 2e$ matrix with elements in
$k[y]$:
\[
\Omega(y)=\left(
\begin{matrix}
\chi_1 & y\chi_3\\
\chi_4 & \chi_2
\end{matrix}
\right).
\]
Since $r$ is the maximum of $\{\rank \Omega_i\}_{i\in\mathbb Z}$,
there exists a nonzero $r\times r$ minor $\Delta_\ell$ of
$\Omega_\ell$ for some $\ell$.  Let $\Delta(y)$ denote the $r \times
r$ minor of $\Omega(y)$ corresponding to $\Delta_\ell$.  Then
$\Delta(y)$ is a polynomial in $y$.  Since
$\Delta(\alpha^\ell)=\Delta_\ell$ is nonzero, $\Delta(y)$ is a nonzero
polynomial. Note that it has degree at most $c(N)$, and therefore it
has at most $c(N)$ roots in $k$. In conclusion, the $r\times r$ minor
$\Delta_i=\Delta(\alpha^i)$ of $\Omega_i$ is zero for at most $c(N)$
values of $i$.
\end{proof}

\begin{remark}
  Formulated in terms of Tate cohomology, cf.\ \ref{Tate}, the proofs
  of parts (1) and (2) of Proposition \ref{3} actually show the
  following:

\begin{enumerate}[\quad\rm(1)]
\item If $\TTor_j^A(M,N)=0$ for some $j$, then $\TTor_i^A(M,N)\ne 0$ for at most $2c(N)$ values of $i$.
\item If $\TExt^j_A(M,N)=0$ for some $j$, then $\TExt^i_A(M,N)\ne 0$ for at most $2c(N)$ values of $i$. 
\end{enumerate}
Similar statements can be given for the ring $B$ and the module $L$.
\end{remark}

\section{Classes of Gorenstein Rings}
\label{new rings}

In this section we discuss homologically defined classes of local 
Gorenstein rings, introduced in \cite{AB} and \cite{HJ}, and show
using the examples in the previous section that these classes of local 
rings lie properly between the class
of local complete intersections and the class of local Gorenstein rings.

Throughout this section $R$ is a local ring with maximal ideal $\fm$.  
Let $\ci$ denote the
condition that $R$ is a local complete intersection, and $\gor$ the
condition that $R$ is a local Gorenstein ring.  We further consider
the following properties (cf. \cite[6.3]{AB}):
\begin{align*}
&\te&&\Tor^R_i(M,N)=0\text{\it\ for all $i\gg 0$ implies\ }
\Ext_R^i(M,N)=0\text{\it\ for all $i\gg 0$},\\
&\et&&\Ext_R^i(M,N)=0\text{\it\ for all $i\gg 0$ implies\ }
\Tor^R_i(M,N)=0\text{\it\ for all $i\gg 0$},\\
&\ee&&\Ext^i_R(M,N)=0\text{\it\ for all $i\gg 0$ implies\ }
\Ext_R^i(N,M)=0\text{\it\ for all $i\gg 0$},
\end{align*}
where $M$ and $N$ range over all finitely generated $R$-modules.  

\begin{chunk}\label{teeegor} 
Note that by taking $N=R$, the property $\te$ implies
$R$ is Gorenstein; by taking $M=R$, the property $\ee$ implies $R$ 
is Gorenstein.
\end{chunk}

\begin{chunk}
\label{AB}
Avramov and Buchweitz prove in \cite[6.1]{AB} that if $R$ is a local
complete intersection, then it satisfies both $\et$ and $\te$.  This gives
the implication (1) in the following diagram, reproduced from \cite[6.3]{AB};
the remaining implications are clear.
$$
\xymatrixrowsep{2pc}
\xymatrixcolsep{2pc}
\xymatrix{
{}&{}&\ee\ar@{<=}[dl]\ar@{=>}[dr]&{}\\
\ci\ar@{=>}^(.35){(1)}[r]&\et{\ \&\ }\te
\ar@{=>}[dr]_{(2)}&{}&\gor\\
{}&{}&\te\ar@{=>}[ur]&{}
}
$$
\end{chunk}

In \cite[6.4]{AB} the question is raised whether
any of these implications can be reversed.
We first note that $(2)$ is reversible:

\begin{proposition}
\label{(3)}
The following statements are equivalent. 
\begin{enumerate}[\quad\rm(1)]
\item $R$ satisfies $\te$. 
\item $R$ is Gorenstein and satisfies $\et$.
\end{enumerate}
\end{proposition}

\begin{proof} 
  We give only the proof of $(1)\implies (2)$. The reverse implication
  can be proved similarly.
  
Assume that $R$ satisfies $\te$. By \ref{teeegor}, the ring $R$ is
Gorenstein. By taking syzygies, it suffices to prove that $\et$ holds
for maximal Cohen-Macaulay $R$-modules $M$ and $N$. Note that such modules 
are in particular reflexive, that is they are isomorphic to their double 
dual.
  
Let $M$ and $N$ be maximal Cohen-Macaulay $R$-modules such that
$\Ext^i_R(M,N)=0$ for all $i\gg 0$. By \cite[2.1]{HJ} we have then
$\Tor_i^R(M,N^*)=0$ for all $i\gg 0$. Since $\te$ holds, this implies
$\Ext^i_R(M,N^*)=0$ for all $i\gg 0$. Using again \cite[2.1]{HJ} we
conclude that $\Tor_i^R(M,N)=0$ for all $i\gg 0$.
\end{proof}

\medskip

In \cite{HJ}, $R$ is said to be an {\it AB ring\/} whenever it satisfies
the following condition:

\medskip

\ab\quad {\it $R$ is Gorenstein and there exists an integer $n$ 
such that for all pairs $(M,N)$ of finitely generated $R$-modules\/} 
$$
\text{\it $\Ext^i_R(M,N)=0$ for all $i\gg 0$ implies $\Ext^i_R(M,N)=0$ for all
$i> n$.}
$$

\medskip

It is shown in \cite{HJ} that if $R$ is an AB ring, then the integer $n$
above can be taken to be $\dim R$, but not less.  Note that the
condition $\ab$ is the Uniform Auslander Condition $\ac$ from
Section \ref{ap} together with the requirement that $R$ is Gorenstein. 

\medskip

Another property of a local Gorenstein ring $R$ is defined in
\cite{HJ}.  We say that $\Ext_R(M,N)$ has a {\it gap of length $g$\/} if
there exists an $n>0$ such that $\Ext^i_R(M,N)\neq 0$ for $i=n-1$ and
$i=n+g$ and $\Ext^i_R(M,N)=0$ for all $n \leq i \leq n+g-1$.
Set
\[
\text{Ext-gap}(R):=\sup\{g\mid \Ext_R(M,N)\text{ has a gap of length
  $g$ }\}\,,
\] where $M$ and $N$ range over all finitely generated
$R$-modules. Similarly, one can define the notion of Tor-gap. It is
proved in \cite[3.3(2)]{HJ} that over a Gorenstein ring Ext-gap is finite
if and only if Tor-gap is finite.

We define the property {\it finite gap} as follows:

\medskip

\gap\quad $R$ {\it is Gorenstein and\/} Ext-gap$(R)$ {\it is finite\/}.

\medskip 

\begin{chunk}
\label{gap}
The following implications are known to hold:
$$
\xymatrixrowsep{2pc}
\xymatrixcolsep{2pc}
\xymatrix{
\ci\ar@{=>}^(.5){(3)}[r]&\gap
\ar@{=>}[r]^{(4)}&\ab\ar@{=>}[r]^{(5)}&\ee
}
$$
The implication (3) is proved in \cite[1.6]{M}, cf. also \cite[2.3]{J}
for a more precise version. The implication (4) is given by
\cite[3.3(3)]{HJ} and
(5) by \cite[4.1]{HJ}.
\end{chunk}

In \cite{HJ}, it is also proved that the implication (3) above and (1)
in \ref{AB} are not reversible. The details of this are as follows.

Every local Gorenstein ring $R$ (which is not a complete intersection) 
has multiplicity at least $\edim
R-\dim R+2$.  A local Gorenstein ring $R$ is said to have {\it minimal
multiplicity} if its multiplicity is equal to $\edim R-\dim R+2$. When
$R$ is artinian, minimal multiplicity just means $\fm^3=0$. A local
Gorenstein ring $R$ of minimal multiplicity is not a complete
intersection precisely when $\edim R-\dim R\geq 3$.

\begin{chunk}
\label{minmult}
Let $R$ be a local Gorenstein ring of minimal multiplicity.  By
\cite[3.6]{HJ} and \cite[3.2(3)]{HJ}, $R$ satisfies $\gap$, and hence
$\ab$.

By \cite[3.6]{HJ}, if $R$ is a Gorenstein ring of minimal multiplicity, 
but not a complete
intersection, then $R$ satisfies the property $\tv$ introduced in
Section \ref{ap}.  Thus all Gorenstein rings of minimal multiplicity
also satisfy $\te$.
\end{chunk}

The main theorem stated in the introduction and proved in the previous
section, shows that 
\begin{enumerate}[\quad (a)]
\item there exist local Gorenstein rings which are not AB rings.
\item there exist local Gorenstein rings which do not satisfy $\te$.
\end{enumerate}

\begin{chunk}
\label{thediag}
The facts discussed above are summarized in the following refinement
of the diagrams in \ref{AB} and \ref{gap}.
$$
\xymatrixrowsep{2.4pc}
\xymatrixcolsep{0.1pc}
\xymatrix{
& & & & &\qquad\et{\ \& \ }\te\ar@{<=>}[rrrr]\ar@{=>}[dlllll]<-0.7ex>|-{\boldsymbol\vert}& & &
&\te
\ar@{=>}[drr]_{(6)}& & & & & & &\\
\ci\ar@{<=}[drrrrr]<-0.7ex>|-{\boldsymbol\vert}\ar@{=>}[urrrrr]<-0.7ex> & & & & &
& & & & & &\ee
\ar@{=>}[rrrr]& & &
&\gor\ar@{=>}[ullllll]|-{\boldsymbol\vert}\ar@{=>}[dllllll]|-{\boldsymbol\vert}\\
& & & & &\gap\ar@{<=}[ulllll]<-0.7ex>\ar@{=>}[rrrr]^{{(4)}} & & &
&
\ab\ar@{=>}[urr]^{(5)} & & & & & & &}
$$
\end{chunk}

\medskip

In conclusion, the following classes of local rings lie strictly
between the class of the local complete intersections and that of 
local Gorenstein rings:
\begin{enumerate}[\quad (a)]
\item the local rings satisfying $\te$.
\item the AB rings
\item the local rings satisfying $\gap$.
\end{enumerate}

\begin{chunk}
Our examples in the previous section are local Gorenstein rings $R$
with $\fm^4=0$, $\edim R=5$, and $\lambda(R)=12$ which satisfy
neither $\ab$ nor $\te$. By \ref{minmult}, these examples are
minimal with respect to the invariant $\inf\{n \mid \fm^n=0\}$. Other
aspects of the minimality of these examples can be deduced from
\cite{S} as we now describe.

Let $R$ be a Gorenstein local ring.  By \cite[3.4]{S} and
\ref{AB}(1), if $\edim R-\dim R\le 4$, then $R$
satisfies $\ab$. Our examples of Gorenstein rings not satisfying
$\ab$ are thus minimal with respect to embedding dimension.

If $R$ is standard graded with $\lambda(R)<12$, then it follows that
$\fm^3=0$ or $\edim R\le 4$, hence
\ref{minmult} or the above considerations apply. In particular, any
such ring satisfies $\ab$, and so in the standard graded case our
examples of rings not satisfying $\ab$ are minimal with respect to length. 
\end{chunk}

\begin{chunk}
There is a question remaining:
$$
\text{Are there further implications between the properties displayed in 
\ref{thediag}?}
$$
\end{chunk}

\section*{Acknowledgements}

The authors thank Luchezar Avramov for his comments and
suggestions on earlier versions of this paper.   They used the
computer algebra package {\em Macaulay 2\/} \cite{M2}
for testing some preliminary examples.


\begin{thebibliography}{99}

\bibitem{A}
M.~Auslander,
{\em Selected works\/},
I.~Reiten, S.~ O.~ Smal{\o}, and {\O}.~ Solberg, Editors,
AMS, Collected Works, vol. 10, 1999.


\bibitem{Av} L.~L.~Avramov, 
{\em Homological asymptotics of modules over local rings}, 
Commutative Algebra (Berkeley, 1987), MSRI
Publ.  {\bf 15}, Springer, New York 1989; pp. 33--62.


\bibitem{AB}
L.~L.~Avramov, R.-O.~Buchweitz,
{\em Support varieties and cohomology over complete intersections\/}, 
Invent. Math. {\bf 142} (2000), 285--318.

\bibitem{ABS}
L.~L.~Avramov, R.-O.~Buchweitz, L.~M.~\c Sega,
{\em Extensions of a dualizing complex by its ring: commutative
versions of a conjecture of Tachikawa\/}, 
J. Pure Appl. Algebra (to appear).

\bibitem{ABr} M.~Auslander, M.~Bridger,
{\em Stable Module Theory\/},
Memoirs of the A.M.S. {\bf 94} (1969),
American Math. Society, Providence, R.I.

\bibitem{AGP} 
L.~L.~Avramov, V.~N.~Gasharov, I.~V.~Peeva,
{\em Complete intersection dimension\/}, 
Publ. Math. I.H.E.S. {\bf 86} (1997), 67--114. 

\bibitem{AY}
T.~Araya, Y.~Yoshino,
{\em Remarks on a depth formula, a grade inequality and a conjecture of 
Auslander\/},
Comm. Alg. {\bf 26} (1998), 3793--3806.

\bibitem{BCR} 
D.~J.~Benson, J.~F.~Carlson, G.~R.~Robinson, 
{\em On the vanishing of group cohomology\/},
J. Algebra {\bf 131} (1990),  40--73.

%\bibitem{CE} H.~Cartan, S.~Eilenberg,
%{\em Homological Algebra\/},
%Princeton University Press (1956), Princeton, NJ.

\bibitem{F} R.~Froberg, 
{\em Koszul algebras\/}, Advances in commutative ring theory 
(Fez, 1997), Lecture Notes in Pure and Appl. Math., 
{\bf 205}, Dekker, New York 1999; pp. 337--350.

\bibitem{GP} 
V.~Gasharov, I.~Peeva, 
{\em Boundedness versus periodicity over commutative local rings\/},
Trans. Amer. Math. Soc. {\bf 320} (1990),  569--580.

\bibitem{H} D.~Happel,
{\em Homological conjectures in representation
theory of finite-dimensional algebras\/}, preprint, 1990.

(Available at:  http://www.math.ntnu.no/\~{}oyvinso/Nordfjordeid/Program/references.html )

\bibitem{Ha} R.~Hartshorne, 
{\em Residues and duality\/}, 
Lecture Notes in Math., {\bf 20}, Springer, Berlin 1971.
 
\bibitem{He} 
R.~Heitmann,
{\em A counterexample to the rigidity conjecture for rings\/},
Bull. Amer. Math. Soc. {\bf 29} (1993), 94--97.

\bibitem{HJ}
C.~Huneke and D.~A.~Jorgensen,
{\em Symmetry in the vanishing of\/ $\Ext$ over Gorenstein rings\/},
Math. Scand. (to appear)

\bibitem{HSV}
C.~Huneke, L.~M.\ \c Sega, A.~N.~Vraciu,
{\em Vanishing of Ext and Tor over Cohen-Macaulay local rings\/},
preprint, 2003.


\bibitem{J}
D.~A.~Jorgensen, 
{\em Vanishing of {\rm(}co{\rm)}homology over commutative rings\/},
Comm. Algebra {\bf 29} (2001), 1883-1898.
 
\bibitem{J1}
D.~A.~Jorgensen, 
{\em A generalization of the Auslander-Buchsbaum formula\/},
J. Pure Appl. Algebra 144 (1999), 145--155. 


\bibitem{M2}
D.~R.~Grayson and M.~E.~Stillman,
{\em Macaulay 2, a software system for research in algebraic geometry\/},
Available at http://www.math.uiuc.edu/Macaulay2

\bibitem{CM}
C.~Miller,
{\em Complexity of tensor products of modules and a theorem of Huneke-Wiegand\/},
Proc. Amer. Math. Soc.
{\bf 126}
(1998),
53--60.


\bibitem{M}
M.~P.~Murthy, 
{\em Modules over regular local rings},
Ill. J. Math. {\bf 7} (1963), 558-565.
        
\bibitem{S}
L.~M.~\c Sega,
{\em Vanishing of cohomology over Gorenstein rings of small codimension\/},
Proc. Amer. Math. Soc. {\bf 131}, (2003), 2313-2323.

\bibitem{Sch}
C.~Scheja,
{\em \"Uber die Bettizahlen lokaler Ringe\/},
Math. Ann.
{\bf 155}
(1964),
155--172.


\end{thebibliography}
\end{document}